\documentclass{amsart}

\usepackage{amsmath,amsthm,amsfonts,amscd,amssymb}
\usepackage[pdfborder={0 0 0}]{hyperref}
\usepackage{enumerate,verbatim}
\usepackage{enumitem}

\newcommand{\ZZ}{\mathbb{Z}}

\newtheorem{thm}{Theorem}[section]
\newtheorem{cor}[thm]{Corollary}
\newtheorem{lem}[thm]{Lemma}

\newtheorem{quest}{Question}

\theoremstyle{definition}

\theoremstyle{remark}
\newtheorem{rem}{Remark}[section]

\begin{document}

\title{On the geometry of cyclic lattices}
\author{Lenny Fukshansky}
\author{Xun Sun}\thanks{The first author was partially supported by NSA Young Investigator Grant \#1210223 and Simons Foundation grants \#208969, 279155.}

\address{Department of Mathematics, 850 Columbia Avenue, Claremont McKenna College, Claremont, CA 91711}
\email{lenny@cmc.edu}
\address{School of Mathematical Sciences, Claremont Graduate University, Claremont, CA 91711}
\email{foxfur\_32@hotmail.com}

\subjclass[2010]{Primary: 11H06, 11H55; Secondary: 68Q17}
\keywords{cyclic lattices, well-rounded lattices, shortest vector problem}

\begin{abstract}
Cyclic lattices are sublattices of $\mathbb Z^N$ that are preserved under the rotational shift operator. Cyclic lattices were introduced by D.~Micciancio in~\cite{mic1} and their properties were studied in the recent years by several authors due to their importance in cryptography. In particular, Peikert and Rosen~\cite{peikert} showed that on cyclic lattices in prime dimensions, the shortest independent vectors problem SIVP reduces to the shortest vector problem SVP with a particularly small loss in approximation factor, as compared to general lattices. In this paper, we further investigate geometric properties of cyclic lattices. Our main result is a counting estimate for the number of well-rounded cyclic lattices, indicating that well-rounded lattices are more common among cyclic lattices than generically. We also show that SVP is equivalent to SIVP on a positive proportion of Minkowskian well-rounded cyclic lattices in every dimension. As an example, we demonstrate an explicit construction of a family of such lattices on which this equivalence holds. To conclude, we introduce a class of sublattices of $\mathbb Z^N$ closed under the action of subgroups of the permutation group $S_N$, which are a natural generalization of cyclic lattices, and show that our results extend to all such lattices closed under the action of any $N$-cycle.
\end{abstract}

\maketitle

\def\A{{\mathcal A}}
\def\AA{{\mathfrak A}}
\def\B{{\mathcal B}}
\def\C{{\mathcal C}}
\def\D{{\mathcal D}}
\def\EE{{\mathfrak E}}
\def\F{{\mathcal F}}
\def\G{{\mathcal G}}
\def\x{{\mathcal H}}
\def\I{{\mathcal I}}
\def\II{{\mathfrak I}}
\def\J{{\mathcal J}}
\def\K{{\mathcal K}}
\def\kk{{\mathfrak K}}
\def\L{{\mathcal L}}
\def\LL{{\mathfrak L}}
\def\M{{\mathcal M}}
\def\mm{{\mathfrak m}}
\def\MM{{\mathfrak M}}
\def\N{{\mathcal N}}
\def\O{{\mathcal O}}
\def\OO{{\mathfrak O}}
\def\PP{{\mathfrak P}}
\def\R{{\mathcal R}}
\def\PNR{{\mathcal P_N(\real)}}
\def\PMNR{{\mathcal P^M_N(\real)}}
\def\PdNR{{\mathcal P^d_N(\real)}}
\def\s{{\mathcal S}}
\def\V{{\mathcal V}}
\def\X{{\mathcal X}}
\def\Y{{\mathcal Y}}
\def\Z{{\mathcal Z}}
\def\H{{\mathcal H}}
\def\cee{{\mathbb C}}
\def\Nn{{\mathbb N}}
\def\pee{{\mathbb P}}
\def\que{{\mathbb Q}}
\def\QQ{{\mathbb Q}}
\def\real{{\mathbb R}}
\def\RR{{\mathbb R}}
\def\zed{{\mathbb Z}}
\def\ZZ{{\mathbb Z}}
\def\aaa{{\mathbb A}}
\def\ff{{\mathbb F}}
\def\HDelta{{\it \Delta}}
\def\kk{{\mathfrak K}}
\def\qbar{{\overline{\mathbb Q}}}
\def\kbar{{\overline{K}}}
\def\ybar{{\overline{Y}}}
\def\kkbar{{\overline{\mathfrak K}}}
\def\ubar{{\overline{U}}}
\def\eps{{\varepsilon}}
\def\beps{{\boldsymbol \varepsilon}}
\def\ahat{{\hat \alpha}}
\def\bhat{{\hat \beta}}
\def\gt{{\tilde \gamma}}
\def\h{{\tfrac12}}
\def\be{{\boldsymbol e}}
\def\bei{{\boldsymbol e_i}}
\def\bc{{\boldsymbol c}}
\def\bm{{\boldsymbol m}}
\def\bk{{\boldsymbol k}}
\def\bi{{\boldsymbol i}}
\def\bl{{\boldsymbol l}}
\def\bq{{\boldsymbol q}}
\def\bu{{\boldsymbol u}}
\def\bt{{\boldsymbol t}}
\def\bs{{\boldsymbol s}}
\def\bv{{\boldsymbol v}}
\def\bw{{\boldsymbol w}}
\def\bx{{\boldsymbol x}}
\def\bX{{\boldsymbol X}}
\def\bz{{\boldsymbol z}}
\def\bwy{{\boldsymbol y}}
\def\bY{{\boldsymbol Y}}
\def\bL{{\boldsymbol L}}
\def\ba{{\boldsymbol a}}
\def\bb{{\boldsymbol b}}
\def\bet{{\boldsymbol\eta}}
\def\bxi{{\boldsymbol\xi}}
\def\bo{{\boldsymbol 0}}
\def\bone{{\boldsymbol 1}}
\def\bol{{\boldsymbol 1}_L}
\def\ep{\varepsilon}
\def\p{\boldsymbol\varphi}
\def\q{\boldsymbol\psi}
\def\rank{\operatorname{rank}}
\def\aut{\operatorname{Aut}}
\def\lcm{\operatorname{lcm}}
\def\sgn{\operatorname{sgn}}
\def\spn{\operatorname{span}}
\def\md{\operatorname{mod}}
\def\Norm{\operatorname{Norm}}
\def\dim{\operatorname{dim}}
\def\det{\operatorname{det}}
\def\Vol{\operatorname{Vol}}
\def\rk{\operatorname{rk}}
\def\ord{\operatorname{ord}}
\def\ker{\operatorname{ker}}
\def\div{\operatorname{div}}
\def\Gal{\operatorname{Gal}}
\def\GL{\operatorname{GL}}
\def\SNR{\operatorname{SNR}}
\def\WR{\operatorname{WR}}
\def\IWR{\operatorname{IWR}}
\def\scg{\operatorname{\left< \Gamma \right>}}
\def\swrh{\operatorname{Sim_{WR}(\Lambda_h)}}
\def\ch{\operatorname{C_h}}
\def\cht{\operatorname{C_h(\theta)}}
\def\scgt{\operatorname{\left< \Gamma_{\theta} \right>}}
\def\scgmn{\operatorname{\left< \Gamma_{m,n} \right>}}
\def\gat{\operatorname{\Omega_{\theta}}}
\def\mn{\operatorname{mn}}
\def\disc{\operatorname{disc}}
\def\rot{\operatorname{rot}}
\def\Prob{\operatorname{Prob}}
\def\co{\operatorname{co}}
\def\Ker{\operatorname{Ker}}

\section{Introduction}
\label{intro}

Define the rotational shift operator on $\real^N$, $N \geq 2$, by
$$\rot(x_1,x_2,\dots,x_{N-1},x_N) = (x_N, x_1,x_2,\dots,x_{N-1})$$
for every $\bx = (x_1,x_2,\dots,x_{N-1},x_N) \in \real^N$. We will write $\rot^k$ for iterated application of $\rot$ $k$ times for each $k \in \zed_{>0}$ (then $\rot^0$ is just the identity map, and $\rot^k = \rot^{N+k}$). It is also easy to see that $\rot$ (and hence each iteration $\rot^k$) is a linear operator. A sublattice $\Gamma$ of $\zed^N$ is called {\it cyclic} if $\rot(\Gamma) = \Gamma$, i.e. if for every $\bx \in \Gamma$, $\rot(\bx) \in \Gamma$. Clearly, $\zed^N$ itself is a cyclic lattice. In fact,  cyclic lattices come from ideals in the quotient polynomial ring $\zed[x]/(x^N -1)$. Let $p(x) \in \zed[x]/(x^N -1)$, then $p(x) = \sum_{n=0}^{N-1} a_n x^n$ for some $a_0,\dots,a_{N-1} \in \zed$. Define a $\zed$-module isomorphism $\rho : \zed[x]/(x^N -1) \to \zed^N$ given by
$$\rho(p(x)) = (a_0,\dots,a_{N-1}) \in \zed^N,$$
then for any ideal $I \subseteq \zed[x]/(x^N -1)$, $\Gamma_I := \rho(I)$ is a sublattice of $\zed^N$. Notice that for every $p(x) = \sum_{n=0}^{N-1} a_n x^n \in I$, 
$$xp(x) = a_{N-1} + a_0x + a_1x^2 + \dots + a_{N-2}x^{N-1} \in I,$$
and so
$$\rho(xp(x)) = (a_{N-1},a_0,a_1,\dots,a_{N-2}) = \rot(\rho(p(x))) \in \Gamma_I,$$
and for any $(a_0,\dots,a_{N-1}) \in \Gamma_I$,
$$\rot(a_0,\dots,a_{N-1}) = \rho \left( x \sum_{n=0}^{N-1} a_n x^n \right) \in \Gamma_I,$$
since $x \sum_{n=0}^{N-1} a_n x^n \in I$. In other words, $\Gamma \subseteq \zed^N$ is a cyclic lattice if and only if $\Gamma = \Gamma_I$ for some ideal $I \subseteq \zed[x]/(x^N -1)$. Cyclic lattices were introduced by D.~Micciancio in~\cite{mic1} and~\cite{mic} in the context of cryptographic algorithms and were further studied in~\cite{lub1}, \cite{peikert}, among other sources. In fact, cyclic lattices are used in the well known NTRU cryptosystem~\cite{ntru}, \cite{ntru1} (also see, for instance~\cite{schneider} and \cite{stehle} for some details) and are further discussed in the context of post-quantum cryptography~\cite{post-quantum}.

On the other hand, given a lattice $\Gamma \subset \real^N$ of rank~$r$, we define its successive minima by
$$\lambda_i = \lambda_i(\Gamma) := \inf \{ \lambda \in \real_{>0} : \Gamma \cap \lambda B_N \text{ contains } i \text{ linearly independent vectors} \},$$
where $B_N$ is a unit ball centered at the origin in $\real^N$, and so
$$0 < \lambda_1 \leq \dots \leq \lambda_r.$$
Let us write $\|\ \|$ for the usual Euclidean norm on~$\real^N$. There exists a collection of linearly independent vectors $\bx_1,\dots,\bx_r$ in $\Gamma$ such that $\|\bx_i\| = \lambda_i$ for each $1 \leq i \leq r$; we will refer to them as vectors {\it corresponding} to successive minima. When $r \leq 4$, there exists a basis for $\Gamma$ consisting of vectors corresponding to successive minima, which is a {\it Minkowski reduced basis} for~$\Gamma$; this is not necessarily true for $r \geq 5$ (see for instance \cite{pohst}), but there are many lattices in higher dimensions as well for which it is true; following J. Martinet, we call such lattices {\it Minkowskian}. Notice also that~$\lambda_1$ is the minimal norm of nonzero vectors in~$\Gamma$ and define the {\it set of minimal vectors}
$$S(\Gamma) = \left\{ \bx \in \Gamma : \|\bx\| = \lambda_1 \right\}.$$
The lattice~$\Gamma$ is called {\it well-rounded} (abbreviated WR) if $\lambda_1 = \dots = \lambda_r$, which is equivalent to saying that $S(\Gamma)$ spans a subspace of~$\real^N$ of dimension~$r$. A strictly stronger condition in general is: $\Gamma = \spn_{\zed} S(\Gamma)$; we will refer to it by saying that $\Gamma$ is $\WR'$. WR lattices are important in discrete optimization, in particular in the investigation of sphere packing, sphere covering, and kissing number problems (see \cite{martinet}), as well as in coding theory (see \cite{esm}). Properties of WR lattices have also been investigated in \cite{mcmullen} in connection with Minkowski's conjecture and in \cite{lf:robins} in connection with the linear Diophantine problem of Frobenius.

Let $\C_N$ be the set of full-rank cyclic sublattices of~$\zed^N$. In this paper we discuss some geometric properties of lattices from~$\C_N$, in particular establishing the following counting estimate on the number of well-rounded cyclic lattices.

\begin{thm} \label{q_cyclic} Let $R \in \real_{> 0}$, then there exists a constant $\alpha_N > 0$ depending only on dimension $N$ such that
\begin{equation}
\label{ratio}
\# \left\{ \Gamma \in \C_N: \lambda_N(\Gamma) \leq R,\ \Gamma \text{ is } \WR'  \right\} \geq \alpha_N R^N
\end{equation}
as $R \to \infty$.
\end{thm}

\begin{rem} \label{mink} By Minkowski Successive Minima Theorem (see, for instance Theorem~2.6.8 on p.~50 of~\cite{martinet}),
$$\det(\Gamma) \gg \ll_N \lambda_N(\Gamma)^N.$$
Hence
$$\# \left\{ \Gamma \in \C_N: \lambda_N(\Gamma) \leq R \right\} \gg \ll_N \# \left\{ \Gamma \in \C_N: \det(\Gamma) \leq R^{\frac{1}{N}} \right\},$$
and analogously for subsets of $\C_N$ consisting of $\WR$ or $\WR'$ lattices.
\end{rem}

\noindent
When $N=2$ a direct argument can be applied to obtain a more explicit bound.

\begin{thm} \label{cyclic_low} Let $R \in \real_{> 0}$, then
\begin{eqnarray}
\label{ratio_2}
0.200650... \times R^2 - 3.035275... \times R & \leq & \# \left\{ \Gamma \in \C_2 : \lambda_2(\Gamma) \leq R,\ \Gamma \text{ is } \WR' \right\} \nonumber \\
& \leq & 0.267638... \times R^2 + 1.673031... \times R.
\end{eqnarray}
\end{thm}

\begin{rem} \label{ideal_comp} The estimate of Theorems~\ref{q_cyclic} and~\ref{cyclic_low} is of the same order of magnitude as the number of {\it all} (not only WR) ideal lattices from polynomial rings $\zed[x]/f(x)$ for irreducible polynomials $f(x)$ under the same map $\rho$ as above (see~\cite{lind}). On the other hand, the number of all cyclic lattices with successive minima $\leq R$ grows like $O(R^N (\log R)^{d(N)-1})$ as $R \to \infty$, where $d(N)$ is the number of divisors of $N$: this is a special case of an estimate of the number of ideal lattices in a forthcoming paper by S. K\"uhnlein and the first author. 
\end{rem}

Lattice-based cryptographic algorithms heavily rely on the fact that the problem of finding $\lambda_1(\Gamma)$, given an arbitrary basis matrix for $\Gamma$, is NP-hard. For most lattices, the problem of finding all successive minima is strictly harder, however if the lattice is WR then the two problems are the same. On the other hand, the set of WR lattices has measure zero in the space of all lattices in a given dimension~$N$. The advantage of using cyclic lattices is that many of them can be constructed from a single vector (using its rotations), and hence the size of the input for a basis matrix of the lattice reduces from~$N^2$ to~$N$. While it is not clear whether the problem of finding $\lambda_1(\Gamma)$ still remains NP-hard, there are reasons to expect that for many cyclic lattices this problem is the same as that of finding all successive minima, i.e. many cyclic lattices are WR. In particular, in~\cite{peikert} the authors proved that in {\it prime} dimensions $N$, the shortest independent vectors problem SIVP on cyclic lattices reduces to (a slight variant of) the shortest vector problem SVP by a polynomial-time algorithm with only a factor of 2 loss in approximation factor (compare to the factor of $\sqrt{N}$ loss on general lattices; see Figure 1 on p.~140 of~\cite{mic_book}). As a corollary of our proof of Theorem~\ref{q_cyclic}, we show that SVP and SIVP are equivalent on a positive proportion of Minkowskian well-rounded cyclic lattices in every dimension $N$ and exhibit a construction of a family of such lattices for which this equivalence holds. These results are given by Lemma~\ref{WR_4}, Remark~\ref{SVP_SIVP} and Corollary~\ref{WR_equiv}.

The paper is organized as follows. In Section~\ref{cyclic-rank} we establish some preliminary results on distribution properties of cyclic lattices. In Section~\ref{gen_cycl} we give a lower bound on the number of $\WR'$ cyclic lattices with bounded successive minima, proving Theorem~\ref{q_cyclic}. Among WR cyclic lattices spanned by their shortest vectors, we specifically focus on those that are in fact spanned by rotations of a single shortest vector: for many such lattices all rotations of any shortest vector are linearly independent, and hence SIVP on these lattices is solved by taking a solution to SVP and all of its rotations. We prove Theorem~\ref{cyclic_low} in Section~\ref{cyclic-2}. Here we follow the tactic of Section~\ref{gen_cycl}, but make the estimates more precise in dimension~2. 

In Section~\ref{perm} we extend our results to a more general class of lattices. Specifically, let $S_N$ be the group of permutations on $N \geq 2$ elements. We can define an action of $S_N$ on $\real^N$ by
\begin{equation}
\label{action}
\tau \bx = \begin{pmatrix} x_{\tau(1)} \\ \vdots \\ x_{\tau(N)} \end{pmatrix}
\end{equation}
for each $\tau \in S_N$ and $\bx = (x_1,\dots,x_N)^t \in \real^N$. We say that a lattice $\Lambda \subset \real^N$ is {\it $\tau$-invariant} (or invariant under $\tau$) for a fixed $\tau \in S_N$ if $\tau \Lambda = \Lambda$. In particular, cyclic lattices are precisely the full-rank sublattices of $\zed^N$ invariant under the $N$-cycle $(1\ 2\dots N)$. The following statement about lattices invariant under arbitrary $N$-cycles follows from our Theorem~\ref{q_cyclic}.

\begin{cor} \label{N-cycle} Let $N \geq 2$, let $\tau \in S_N$ be an $N$-cycle, and let $\C_N(\tau)$ be the set of all $\tau$-invariant full-rank sublattices of $\zed^N$. Then
\begin{equation}
\label{ratio-1}
\# \left\{ \Gamma \in \C_N(\tau): \lambda_N(\Gamma) \leq R,\ \Gamma \text{ is } \WR'  \right\} \geq \alpha_N R^N,
\end{equation}
as $R \to \infty$, for the same value of $\alpha_N$ as in~\eqref{ratio}.
\end{cor}

We prove Corollary~\ref{N-cycle} in Section~\ref{perm} and conclude with some further questions about more general permutation invariant lattices. We are now ready to proceed.
\bigskip

\section{Basic properties of cyclic lattices}
\label{cyclic-rank}

Let $\G_N$ be the set of full-rank cyclic sublattices of~$\zed^N$ spanned by vectors corresponding to their successive minima (when $N \leq 4$, $\G_N = \C_N$). In this section we start out by looking at the cyclic lattices generated by rotations of a single vector. Notice that for every $\ba \in \zed^N$, $\|\ba\| = \|\rot(\ba)\|$, therefore if $\Gamma \subseteq \zed^N$ is a cyclic lattice and $\ba \in S(\Gamma)$, then $\rot^n(\ba) \in S(\Gamma)$ for every $1 \leq n \leq N-1$ (clearly $\rot^N(\ba) = \ba$). Therefore cyclic lattices have large sets of minimal vectors, and so it is natural to expect that they are WR fairly often. In fact, it is clear that if $\ba \in S(\Gamma)$ and $\ba,\rot(\ba),\dots,\rot^{N-1}(\ba)$ are linearly independent, then $\Gamma$ is WR. To state our first observation in this direction, we need some more notation. 

Let $\ba = (a_0,\dots,a_{N-1})^t \in \real^N$, and define $\ba(x) = \sum_{n=0}^{N-1} a_n x^n$ to be the polynomial of degree $ \leq N-1$ in $x$ whose coefficient vector is $\ba$. Let also
$$M(\ba) = (\ba\ \rot(\ba)\ \dots\ \rot^{N-1}(\ba))$$
be an $N \times N$ matrix. Consider the lattice
$$\Lambda(\ba) = \spn_{\zed} \left\{ \ba, \rot(\ba), \dots, \rot^{N-1}(\ba) \right\} = M(\ba) \zed^N,$$
and define the {\it cyclic order} of $\ba$, denoted $\co(\ba)$, to be the rank of $\Lambda(\ba)$. This means that precisely $\co(\ba)$ of the vectors $\ba, \rot(\ba), \dots, \rot^{N-1}(\ba)$ are linearly independent, and so $M(\ba)$ is a matrix of rank $\co(\ba)$. While not every $\Lambda(\ba)$ is necessarily generated by the vectors corresponding to its successive minima, lattices of the form~$\Lambda(\ba)$ for $\ba \in \zed^N$ are very common among cyclic lattices.

\begin{lem} \label{cycl} The vectors $\ba,\rot(\ba),\dots,\rot^{N-1}(\ba) \in \zed^N$ are linearly independent if and only if the polynomial $\ba(x)$ does not have any common factors with $x^N-1$.
\end{lem}

\proof
In this case $M(\ba)$ is an $N \times N$ circulant matrix corresponding to a vector $\ba \in \zed^N$. It is a well-known fact (see for instance~\cite{circ}) that
$$\det(M(\ba)) = \prod_{n=0}^{N-1} \ba(\omega_j),$$
where $\omega_j = e^{\frac{2\pi ij}{N}}$ is an $N$-th root of unity. Hence $\det(M(\ba)) = 0$ if and only if $\ba(\omega_j) = 0$ for some $0 \leq j \leq N-1$, which happens if and only if $\ba(x)$ is divisible by the minimal polynomial of $\omega_j$ -- that is, by some cyclotomic polynomial dividing $x^N-1$.
\endproof

\begin{rem} \label{pkrt} An immediate consequence of Lemma~\ref{cycl} is that when~$N$ is prime, the vectors $\ba,\rot(\ba),\dots,\rot^{N-1}(\ba) \in \zed^N$ are linearly independent if and only if $\ba(x)$ is not a multiple of $x-1$ or~$\sum_{n=0}^{N-1} x^n$. See Section~2 of~\cite{peikert} for further results of this kind. 
\end{rem}

Let
$$C_R^N = \{ \bx \in \real^N\ :\ |\bx| := \max \{|x_1|,\dots,|x_N|\} \leq R \}$$
for every $R \in \real_{>0}$, i.e., $C_R^N$ is a cube of side-length $2R$ centered at the origin in~$\real^N$. Recall that $d$-th cyclotomic polynomial $\Phi_d(x)$ divides $x^N-1$ if and only if $d$ is a divisor of $N$. For each divisor $d$ of $N$, define the {\it $d$-th cyclotomic subspace} to be
\begin{equation}
\label{H_phi}
H_{\Phi_d} = \left\{ \ba \in \real^N : \Phi_d(x) \text{ divides } \ba(x) \text{ in } \real[x] \right\}.
\end{equation}
By Lemmas 2.3 and 2.4 of \cite{peikert}, $H_{\Phi_d}$ is a subspace of $\real^N$ of dimension
$$\dim_{\real} (H_{\Phi_d}) = N-\deg(\Phi_d) = N-\varphi(d),$$
where $\varphi$ is Euler's $\varphi$-function. Then $\Lambda_{\Phi_d} := H_{\Phi_d} \cap \zed^N$ is a sublattice of $\zed^N$ of rank $N-\varphi(d)$. Therefore
\begin{eqnarray}
\label{est}
\left| C_R^N \cap \left( \zed^N \setminus \bigcup_{d \mid N} \Lambda_{\Phi_d} \right) \right| & = & \left| C_R^N \cap \zed^N \right| - \sum_{d \mid N} \left| C_R^N \cap \Lambda_{\Phi_d} \right| \nonumber \\
& \geq & \left| C_R^N \cap \zed^N \right| - \sum_{d \mid N} \left| C_R^{N-\varphi(d)} \cap \zed^{N-\varphi(d)} \right| \nonumber \\
& \geq & \left| C_R^N \cap \zed^N \right| - \left| C_R^{N-1} \cap \zed^{N-1} \right| \sum_{d \mid N} \varphi(d) \nonumber \\
& = & (2R+1)^N - N(2R+1)^{N-1} \nonumber \\
& = & (2R+1-N)(2R+1)^{N-1}.
\end{eqnarray}

\noindent
The lattice $\Lambda(\ba) \subseteq \zed^N$ has rank $N$ if and only if the vectors $\ba,\rot(\ba),\dots,\rot^{N-1}(\ba)$ are linearly independent, which happens if and only if the polynomial $\ba(x)$ is not divisible by any cyclotomic polynomial $\Phi_d(x)$ for any $d \mid N$, by Lemma~\ref{cycl}. How often does this happen?

\begin{lem} \label{cycl_wr_prob} Let $R > \frac{N-1}{2}$, then
\begin{equation}
\label{prb}
\Prob_{\infty,R} \left( \rk(\Lambda(\ba)) = N \right) \geq 1 - \frac{N}{2R+1},
\end{equation}
where probability $\Prob_{\infty,R}( \cdot)$ is with respect to the uniform distribution among all points $\ba$ in the set $C_R^N \cap \zed^N$.
\end{lem}

\proof
By Lemma~\ref{cycl},
$$\Prob_{\infty,R} \left( \rk(\Lambda(\ba)) = N \right) = \frac{\left| C_R^N \cap \left( \zed^N \setminus \bigcup_{d \mid N} \Lambda_{\Phi_d} \right) \right|}{\left| C_R^N \cap \zed^N \right|},$$
and the statement of the lemma follows by~\eqref{est} combined with the observation that $\left| C_R^N \cap \zed^N \right| = (2R+1)^N$.
\endproof
\bigskip

\section{General cyclic lattices}
\label{gen_cycl}

The main goal of this section is to prove Theorem~\ref{q_cyclic}. Recall that $\C_N$ is the set of all cyclic full-rank sublattices of~$\zed^N$, while $\G_N \subset \C_N$ is the subset consisting of all lattices in $\C_N$ which are spanned by the vectors corresponding to successive minima. Naturally, every lattice $\Gamma \in \C_N$ has a sublattice $\Gamma_1 \in \G_N$ which is spanned by the vectors corresponding to successive minima of $\Gamma$; it is called a {\it Minkowskian sublattice} of $\Gamma$. While Minkowskian sublattice may not be unique, there can only be finitely many of them, where an upper bound on this number depends only on~$N$. On the other hand, the index $|\Gamma : \Gamma_1|$ of a Minkowskian sublattice is also bounded above by a constant depending only on~$N$, and hence a given lattice in $\G_N$ can be a Minkowskian sublattice for only finitely many lattices in~$\C_N$ (see~\cite{martinet1} and subsequent works of J. Martinet and his co-authors for more information on the index of Minkowskian sublattices). This means that the numbers of WR lattices in $\C_N$ and $\G_N$ have the same asymptotic order. Here we will construct large families of WR lattices in~$\G_N$.

For a subspace $V \subseteq \real^N$ which is closed under the rotational shift operator, define the set
\begin{equation}
\label{D_def}
\D^V_N = \left\{ \ba \in V : \co(\ba) = \dim_{\real}(V),\ \ba \in S(\Lambda(\ba)),\ \Lambda(\ba) \text{ spanned by } S(\Lambda(\ba)) \right\},
\end{equation}
and let us write $\D_N$ for $\D_N^{\real^N}$.

\begin{lem} \label{D_param} A lattice $\Lambda(\ba) \subset V \subseteq \real^N$ is of rank $= \dim_{\real}(V)$ with $\ba \in S(\Lambda(\ba))$ if and only if $\ba \in \D^V_N$. Moreover, $\Lambda(\ba) = \Lambda(\bb)$ for only finitely many~$\bb \in \D^V_N$ with an upper bound on their number, call it~$\beta(V)$, depending only on the dimension of~$V$; we will write $\beta_N$ for $\beta(\real^N)$.
\end{lem}

\proof
The first assertion is clear from the definition of $\D^V_N$. The second assertion follows from a well known fact in the reduction theory of positive definite quadratic forms (see, for instance, Theorems~1.1-1.2 in Chapter~12 of~\cite{cassels_book}).
\endproof

For each $R \in \real_{>0}$, let $B^V_N(R)$ be a ball of radius $R$ centered at the origin in $V$, and let
$$\D^V_N(R) = \left\{ \ba \in \D^V_N : \|\ba\| \leq R \right\} = \D^V_N \cap B^V_N(R).$$
It is easy to notice that $\ba \in  \D^V_N$ if and only if $R\ba \in \D^V_N$, and hence $\D^V_N(R) = R\D^V_N(1)$ is a homogeneously expanding domain. Moreover, $\D_N^V(R)$ is a symmetric bounded star body, and hence is Jordan-measurable. We write $\D_N(R)$ for $\D_N \cap B_N(R)$, where $B_N(R)$ is a ball of radius $R$ centered at the origin in $\real^N$.

Given a vector $\ba \in \real^N$ with $\co(\ba)=k$, let $\ba_1,\dots,\ba_k$ be some fixed ordering of the vectors $\ba,\rot(\ba),\dots,\rot^{k-1}(\ba)$. Define the angle sequence $\{ \theta_1,\dots,\theta_{k-1} \}$ of this ordering as follows: for each $1 \leq i \leq k-1$, let $\theta_i$ be the angle between $\ba_{i+1}$ and the subspace spanned by $\ba_1,\dots,\ba_i$.

\begin{lem} \label{R-N} Let $V \subseteq \real^N$ be an $L$-dimensional subspace closed under the rotational shift operator. Assume that $V$ contains a vector $\ba$ with $\co(\ba)=L$ such that some ordering of its $L$ linearly independent rotations has the corresponding angle sequence satisfying the condition
\begin{equation}
\label{angles}
\pi/3 + \eps \leq \theta_i \leq 2\pi/3 - \eps
\end{equation}
for each $1 \leq i \leq k-1$, for some $\eps > 0$. Then $\Vol_L(\D^V_N(R)) = O(R^L)$,
where the constant in the $O$-notation depends on $V$, $L$, and $N$.
\end{lem}

\proof
Let $\ba_1,\dots,\ba_L$ be the ordering of $L$ linearly independent rotations of $\ba$ with the corresponding angle sequence as in~\eqref{angles}. Notice that $\|\ba_1\| = \dots = \|\ba_L\| = \|\ba\|$, and so Theorem~1 of~\cite{near:ort} guarantees that $\ba_1,\dots,\ba_L$ are minimal vectors in $\Lambda(\ba)$, hence $\ba \in \D^V_N$.

Let $\delta > 0$ and let 
$$B(V,\delta) = \left\{ \bx \in V : \|\bx\| \leq \delta \right\}$$
be the closed ball of radius $\delta$ centered at the origin in $V$. Let $\bt \in B(V,\delta)$ and $\ba' = \ba+\bt$. Let $\ba'_1,\dots,\ba'_L$ be the rotations of $\ba'$ corresponding to the rotations $\ba_1,\dots,\ba_L$ of $\ba$. There exists a $\delta > 0$, depending on $\eps$, small enough so that for every $\bt \in B(V,\delta)$ the angle sequence $\{ \theta'_1,\dots,\theta'_{k-1} \}$ of $\ba'_1,\dots,\ba'_L$ still satisfies~\eqref{angles} with $\eps$ replaced by some $\eps' > 0$. Then, as above, Theorem~1 of~\cite{near:ort} guarantees that $\ba' \in \D^V_N$, i.e., $\ba + B(V,\delta) \subseteq \D^V_N$, and so $\D^V_N$ must have positive $L$-dimensional volume. Since $\D^V_N$ is a homogeneously expanding domain, we must have
$$0 < \Vol_L(\D^V_N(R)) = \Vol_L(R \D^V_N(1)) = O(R^L),$$
which completes the proof of the lemma.
\endproof

\begin{rem} \label{RN_1} We will apply Lemma~\ref{R-N} to~$\real^N$. Notice that the angle sequence of the rotations of the first standard basis vector $\be_1 \in \real^N$ satisfies the assumption of Lemma~\ref{R-N}. Hence $\Vol_N(\D_N(R)) = O(R^N)$ for every $N \geq 2$, by Lemma~\ref{R-N}.
\end{rem}

\begin{rem} \label{RN_2} There is also another way to look at the set $\D^V_N$ with $V$ as in the statement of Lemma~\ref{R-N}. For each $\ba \in V$, all rotations of $\ba$ have to be in $V$, and so $\co(\ba) \leq L$. Let
\begin{equation}
\label{rot-0}
M_V(\ba) = (\ba\ \rot(\ba)\ \dots\ \rot^{L-1}(\ba)),
\end{equation}
and notice that $M_V(\ba) = M(\ba)$ when $V=\real^N$. Define the corresponding $L \times L$ Gram matrix
$$Q_V(\ba) = M_V(\ba)^t M_V(\ba),$$
and let us write $q_{ij}$ for the entires of this matrix, then
$$q_{ij} = q^V_{ij}(\ba) := \rot^{i-1}(\ba) \cdot \rot^{j-1}(\ba).$$
Notice that
\begin{equation}
\label{rot-1}
\rot^{i-1}(\ba) \cdot \rot^{j-1}(\ba) = \rot^i(\ba) \cdot \rot^j(\ba),
\end{equation}
and so all the distinct entries $q_{ij}$ are represented in the first row. Furthermore,
\begin{equation}
\label{rot-2}
\ba \cdot \rot^{i-1}(\ba) = \ba \cdot \rot^{N-i+1}(\ba)
\end{equation}
for each $2 \leq i \leq N-1$, and hence the total number of distinct off-diagonal entries in the matrix~$Q_V(\ba)$ is at most $[N/2]$; all the diagonal entries $q_{ii} = \|\ba\|^2$. Now, $\ba \in \D^V_N$ if and only if $Q_V(\ba)$ is in the corresponding Minkowski reduction domain, which is known to be a convex polyhedral cone in~$\real^{\frac{L(L+1)}{2}}$ with a finite number of facets (see, for instance, Chapter~12 of~\cite{cassels_book} or~\cite{achill_book}), and conditions~\eqref{rot-0}, \eqref{rot-1}, \eqref{rot-2} imply that $Q_V(\ba)$ would have to be in a specific section of this cone. On the other hand, given a Gram matrix $Q$, the basis matrix $M$ such that $Q=M^tM$ is uniquely determined up to an orthogonal transformation.
\end{rem}

\begin{lem} \label{WR_1} Let $R \in \real_{>0}$, and define
\begin{equation}
\label{f-N}
f_N(R) = \# \left\{ \Lambda(\ba) \in \C_N : \|\ba\| = \lambda_1(\Lambda(\ba)) =  \lambda_N(\Lambda(\ba)) \leq R \right\},
\end{equation}
then
\begin{equation}
\label{f-N-bound}
O(R^N) \leq f_N(R) \leq O(R^N),
\end{equation}
where the constants in the $O$-notation depend only on $N$.
\end{lem}

\proof
Let $\beta_N$ be as in Lemma~\ref{D_param}, then
\begin{equation}
\label{f_bound}
\frac{1}{\beta_N} \# \left( \zed^N \cap \D_N(R) \right) \leq f_N(R) \leq \# \left( \zed^N \cap \D_N(R) \right)
\end{equation}
by Lemma~\ref{D_param}. Theorem~2 on p.~128 of~\cite{lang} asserts that
\begin{equation}
\label{b1}
\# \left( \zed^N \cap \D_N(R) \right) = \Vol_N(\D_N(R)) + O(R^{N-1}).
\end{equation}
and so~\eqref{f-N-bound} follows by combining~\eqref{b1} with Lemma~\ref{R-N} and~\eqref{f_bound}.
\endproof

\begin{rem} \label{lipschitz} The boundary of the set $\D_N(R)$ is Lipschitz parameterizable, however that is not important for the application of Theorem~2 on p.~128 of~\cite{lang} in the argument above, since we are only using the main term of the asymptotic formula in our inequalities, and Lemma~\ref{R-N} implies that there exist sets $C_1$, $C_2$ with Lipschitz parameterizable boundaries (in fact, convex sets) such that $RC_1 \subseteq \D_N^V(R) \subseteq RC_2$ for all $R > 0$.
\end{rem}

\proof[Proof of Theorem~\ref{q_cyclic}]
The theorem now follows from the estimates of Lemma~\ref{WR_1}.
\endproof

Now we comment on the connection of our results to the equivalence of SVP and SIVP. Let
$$\R_N = \left\{ \Lambda(\ba) \in \C_N : \|\ba\| = \lambda_1(\Lambda(\ba)) =  \lambda_N(\Lambda(\ba)) \right\},$$
and let $\Gamma \in \R_N$. Suppose that $\bc, \rot(\bc),\dots,\rot^{N-1}(\bc)$ are linearly independent for every $\bc \in S(\Gamma)$, then SIVP is equivalent to SVP on $\Gamma$. In the next lemma we prove that this is true for a positive proportion of lattices in~$\R_N$. Specifically, let
$$\R'_N = \left\{ \Gamma \in \R_N : \co(\bc) = N\ \forall\ \bc \in S(\Gamma) \right\},$$
and define
$$f'_N(R) = \# \left\{ \Gamma \in \R'_N : \lambda_N(\Gamma) \leq R \right\}$$
for any $R \in \real_{> 0}$.

\begin{lem} \label{WR_4} As $R \to \infty$, we have
$$\frac{f'_N(R)}{f_N(R)} \geq O(1),$$
where the constant in $O$-notation depends only on $N$.
\end{lem}

\proof
Let $\Gamma \in \R_N$, and suppose that $\bc \in S(\Gamma)$ is such that $\co(\bc) < N$. Then $\bc \in \Gamma \cap H_{\Phi_d}$ for some $d \mid N$. In other words, $\Gamma \in \R_N \setminus \R'_N$ if and only if
\begin{equation}
\label{pol_co}
S(\Gamma) \cap \left( \bigcup_{d \mid N} H_{\Phi_d} \right) \neq \emptyset.
\end{equation}
Then
$$f'_N(R) \asymp \# \left\{ \ba \in \zed^N \cap \D_N(R) : \Gamma = \Lambda(\ba) \text{ does not satisfy \eqref{pol_co}} \right\},$$
and since~\eqref{pol_co} is given by finitely many polynomial conditions, we have $f'_N(R)  \asymp f_N(R)$.
\endproof

\begin{rem} \label{SVP_SIVP} Lemma~\ref{WR_4} then guarantees that
\begin{equation}
\label{ratio_s}
\frac{\# \left\{ \Gamma \in \R'_N: \lambda_N(\Gamma) \leq R \right\}}{\# \left\{ \Gamma \in \R_N: \lambda_N(\Gamma) \leq R \right\}} \geq O(1) \text{ as } R \to \infty.
\end{equation}
By our observation above, SVP and SIVP are equivalent on $\R'_N$, and so the two problems are equivalent on a positive proportion of cyclic lattices in~$\R_N$.
\end{rem}
\smallskip

In fact, we can use the idea in the proof of Lemma~\ref{R-N} and Remark~\ref{RN_1} to explicitly construct full-rank WR lattices of the form $\Lambda(\ba)$ in~$\real^N$ on which SVP and SIVP are equivalent.

\begin{cor} \label{WR_equiv} Let $k_1,\dots,k_{N-1} \in \zed$ be nonzero integers, $m = \lcm(k_1, \dots, k_{N-1})$, and 
$$\ba = \left( m, \frac{m}{k_1},\dots,\frac{m}{k_{N-1}} \right)^t \in \zed^N.$$
There exists a sufficiently large positive integer $l$, depending only on the dimension~$N$, such that whenever $|k_1|,\dots,|k_{N-1}| \geq l$, the lattice $\Lambda(\ba) \in \R'_N$.
\end{cor}

\proof
Let $l$ be a positive integer, the choice of which is to be specified below, and let the rest of the notation be as in the statement of the corollary. Let $\bb = \frac{1}{m} \ba = \be_1 + \beps$, where
$$\beps = (0, 1/k_1,\dots,1/k_{N-1}).$$
Taking $l$ sufficiently large, we can ensure that the angle sequence of the rotations of the vector~$\bb$ satisfies condition~\eqref{angles} for some $\eps>0$, in which case~$\Lambda(\bb)$ is a lattice of rank~$N$ with minimal norm equal to~$\|\bb\|$ by the same argument as in the proof of Lemma~\ref{R-N} and Remark~\ref{RN_1}. 

We can assume that $l > 10N$ so that $(1- N/l)^2 > 81/100$. We will now show that
\begin{equation}
\label{SLb}
S(\Lambda(\bb)) = \{ \pm \bb, \pm \rot(\bb),\dots, \pm \rot^{N-1}(\bb) \}.
\end{equation}
Indeed, suppose
$$\bc = \sum_{i=1}^{N} \alpha_i \rot^{i-1}(\bb) \in S(\Lambda(\bb)),$$
where $\alpha_1,\dots,\alpha_{N} \in \zed$, not all zero. Let $\alpha = \max_{1 \leq i \leq N} |\alpha_i|$, so for each $1 \leq n \leq N$
$$\left| \alpha_1 + \dots + \alpha_{n-1} + \alpha_{n+1} + \dots + \alpha_N \right| \leq N \alpha.$$
Then $c_n$, the $n$-th coordinate of $\bc$, satisfies the inequalities
$$\max \{ 0, |\alpha_n| - N \alpha/l \} \leq |c_n| \leq |\alpha_n| + N \alpha/l,$$
and so we have
$$\|\bc\|^2 \geq \alpha^2 (1- N/l)^2.$$
Assume first that $\alpha > 1$, then we have
$$\|\bc\|^2 > 2 > 1 + (N-1)/l^2 \geq \|\bb\|^2.$$
Therefore we must have $\alpha=1$. If $\alpha_n = \pm 1$ for only one $n$, then $\bc = \pm \rot^{n-1}(\bb)$. Hence assume there exist $1 \leq j < n \leq N$ such that $\alpha_j, \alpha_n = \pm 1$, then
$$\|\bc\|^2 \geq 2 (1- N/l)^2 > 1 + (N-1)/l^2 = \|\bb\|^2,$$
which establishes~\eqref{SLb}. Then $\Lambda(\ba) = m \Lambda(\bb)$, and hence 
$$S(\Lambda(\ba)) = \{ \pm \ba, \pm \rot(\ba),\dots, \pm \rot^{N-1}(\ba) \},$$
meaning that each vector in $S(\Lambda(\ba))$ has cyclic order $=N$. Thus $\Lambda(\ba) \in \R'_N$.
\endproof

\begin{rem} \label{WR_equiv-1} To summarize, the main idea of Corollary~\ref{WR_equiv} is to pick a rational vector~$\bb$ from a small ball centered at $\be_1$. Then the set of minimal vectors of~$\Lambda(\bb)$ will consist only of $\pm$ rotations of~$\bb$ due to the fact that one coordinate of~$\bb$ strongly dominates others. Hence SVP and SIVP are equivalent on~$\Lambda(\bb)$, and $\Lambda(\bb)$ is similar to some full-rank WR cyclic sublattice of~$\zed^N$ because coordinates of $\bb$ are rational. Since a ball of positive radius centered at $\be_1$ contains infinitely many rational points, infinitely many mutually non-similar lattices with this equivalence property can be constructed this way.
\end{rem}

\bigskip

\section{Cyclic lattices in the plane}
\label{cyclic-2}

In this section we prove Theorem~\ref{cyclic_low}. Recall that every planar cyclic lattice is spanned by vectors corresponding to its successive minima. Furthermore, for a sublattice $\Gamma$ of $\zed^2$, $|S(\Gamma)|=2$ or 4, and $\Gamma$ is WR if and only if $|S(\Gamma)| = 4$. If $\Gamma$ is not WR, then $|S(\Gamma)|=2$ and the vectors corresponding to first and second successive minima are unique (up to $\pm$ sign): this follows, for instance, from the second Theorem and discussion after it on p. 203 of~\cite{donaldson}.

\begin{lem} \label{cycl-2} A lattice $\Gamma \in \C_2$ is WR if and only if either $\Gamma = \Lambda(\ba)$ for some $\ba \in S(\Gamma)$ or $\Gamma =\alpha \begin{pmatrix} 1 & 1 \\ 1 & -1 \end{pmatrix} \zed^2$ for some $\alpha \in \zed_{>0}$. On the other hand, $\Gamma \in \C_2$ is not WR if and only if $\Gamma = \begin{pmatrix} \alpha & \beta \\ \alpha & -\beta \end{pmatrix} \zed^2$ for some distinct positive integers $\alpha, \beta$.
\end{lem}

\proof
If $\Gamma = \Lambda(\ba)$ for some $\ba \in S(\Gamma)$, then $S(\Gamma) = \{ \pm \ba, \pm \rot(\ba)\}$ and the vectors $\ba, \rot(\ba)$ are linearly independent. If  $\Gamma =\alpha \begin{pmatrix} 1 & 1 \\ 1 & -1 \end{pmatrix} \zed^2$ for some $\alpha \in \zed$, then 
$$S(\Gamma) = \left\{ \pm \alpha \begin{pmatrix} 1 \\ 1 \end{pmatrix}, \pm \alpha \begin{pmatrix} 1 \\ -1 \end{pmatrix} \right\}.$$
In both cases, it is clear that $\Gamma$ is WR. 
 
Suppose then that $\Gamma$ is WR, then $|S(\Gamma)| = 4$ and $S(\Gamma)$ contains a basis for $\Gamma$. Let $\ba \in S(\Gamma)$. First assume $\Lambda(\ba)$ has rank 2, then $\ba, \rot(\ba) \in S(\Gamma)$ are linearly independent, and hence form a basis for $\Gamma$. Therefore $\Gamma=\Lambda(\ba)$. Next suppose that $\Lambda(\ba)$ has rank 1, then $\ba = c \rot(\ba)$ for some $c \in \zed$, which easily implies that $a_1=a_2$, and so $\ba = \alpha \begin{pmatrix} 1 \\ 1 \end{pmatrix}$ for some $\alpha \in \zed$. Since $\Gamma$ is WR, there must exist $\bc \in S(\Gamma)$ such that $\bc \neq \pm \ba$. Then $\rot(\bc)$ is also in $S(\Gamma)$, and since $|S(\Gamma)| = 4$, we must have $-\bc = \rot(\bc)$ and $\|\bc\| = \|\ba\|$, meaning that $\bc = \alpha \begin{pmatrix} -1 \\ 1 \end{pmatrix}$. Then $S(\Gamma) = \{ \pm \ba, \pm \bc\}$, and so
$$\Gamma = \alpha \begin{pmatrix} 1 & 1 \\ 1 & -1 \end{pmatrix} \zed^2.$$
This completes the proof of the first statement.

The second statement follows immediately from the observation that~$\real^2$ has precisely two cyclotomic subspaces:
$$H_{\Phi_1} = \spn_{\real} \left\{ \begin{pmatrix} 1 \\ 1 \end{pmatrix} \right\},\ H_{\Phi_2} = \spn_{\real} \left\{ \begin{pmatrix} 1 \\ -1 \end{pmatrix} \right\}.$$
\endproof

For $R \in \real_{>0}$, let $f_2(R)$ be as in~\eqref{f-N} for $N=2$, and define
$$g_2(R) = \# \left\{ \Gamma \in \C_2 : \Gamma \neq \Lambda(\ba)\ \forall\ \ba \in \zed^2,\ \lambda_1(\Gamma) =  \lambda_2(\Gamma) \leq R \right\}.$$
We can now use Lemma~\ref{cycl-2} to estimate the functions~$f_2(R)$ and $g_2(R)$.

\begin{lem} \label{cycl-2-est} Let $R \in \real_{>0}$, then
\begin{equation}
\label{f-2}
0.200650... \times R^2 - 3.742382... \times R \leq f_2(R) \leq 0.267638... \times R^2 + 0.965925... \times R,
\end{equation}
\begin{equation}
\label{g-2}
g_2(R) = \left[ \frac{R}{\sqrt{2}} \right].
\end{equation}
\end{lem}

\proof
First assume $\Gamma = \Lambda(\ba)$ for some $\ba = \begin{pmatrix} a_1 \\ a_2 \end{pmatrix} \in S(\Gamma)$. Notice that we can assume without loss of generality that $|a_1| > |a_2|$. The condition that $\ba,\rot(\ba)$ form a Minkowski reduced basis amounts to satisfying the following condition (see, for instance, Note~1 on p.~257 of~\cite{cassels_book}):
$$a_1^2+a_2^2 \geq 4|a_1a_2|.$$
This means that either
\begin{equation}
\label{C1}
a_1^2+a_2^2-4a_1a_2 \geq 0,\ a_1a_2 \geq 0,
\end{equation}
or
\begin{equation}
\label{C2}
a_1^2+a_2^2+4a_1a_2 \geq 0,\ a_1a_2 < 0.
\end{equation}
First consider the~\eqref{C1} situation, then there are the following two options:
\begin{enumerate}
\item $a_1 \geq [(2+\sqrt{3})a_2] +1 > a_2 \geq 0$,
\item $0 \geq a_2 > [(2+\sqrt{3})a_2] - 1 \geq a_1$.
\end{enumerate}
Notice that $a_1,a_2$ satisfy option (1) if and only if $-a_1,-a_2$ satisfy option (2), hence they correspond to the same lattice $\Lambda(\ba)$. Next consider the~\eqref{C2} situation, then there are the following two options:
\begin{enumerate}[resume]
\item $a_1 \leq -[(2+\sqrt{3})a_2] -1 < 0 < a_2$,
\item $a_1 \geq -[(2+\sqrt{3})a_2] + 1 > 0 > a_2$.
\end{enumerate}
Again, $a_1,a_2$ satisfy option (3) if and only if $-a_1,-a_2$ satisfy option (4), hence they correspond to the same lattice $\Lambda(\ba)$. Notice also that for each pair $a_1,a_2$ satisfying options (1) and (2), there is precisely one pair satisfying options (3) and (4). Hence we will only count vectors $\ba \in \zed^2$ with $\|\ba\| \leq R$ satisfying (1) and multiply this number by 2. Therefore:
\begin{equation}
\label{f_2_R_sum}
f_2(R) = 2 \sum_{a_2=1}^{A(R)} \left( \left[ \sqrt{R^2-a_2^2} \right] - \left[ (2+\sqrt{3})a_2 \right] - 1 \right),
\end{equation}
where
$$A(R) = \left[ \frac{R}{2\sqrt{2+\sqrt{3}}} \right].$$
Using~\eqref{f_2_R_sum}, we now give quick estimates on $f_2(R)$. A higher degree of precision is easily possible here, but we choose in favor of simplicity. Notice that
\begin{eqnarray}
\label{l2-1}
f_2(R) & \geq & 2 R A(R) - 2(3+\sqrt{3}) \sum_{a_2=1}^{A(R)} a_2 - 2A(R) \nonumber \\
& = & 2RA(R) - (3+\sqrt{3}) A(R)^2 - (5+\sqrt{3})A(R) \nonumber \\
& \geq & \frac{\left( 4\sqrt{2+\sqrt{3}} - 3 - \sqrt{3} \right) R^2}{8+4\sqrt{3}} - \frac{\left( 5+\sqrt{3} + 4 \sqrt{2+\sqrt{3}} \right)R}{2\sqrt{2+\sqrt{3}}} \nonumber \\
& = & 0.200650... \times R^2 - 3.742382... \times R.
\end{eqnarray}
On the other hand,
\begin{eqnarray}
\label{l2-2}
f_2(R) & \leq & 2RA(R) - (2+\sqrt{3})A(R)(A(R)+1) \nonumber \\
& \leq & \frac{R^2}{\sqrt{2+\sqrt{3}}} - \frac{R^2}{4} +  \frac{\sqrt{2+\sqrt{3}}\ R}{2} \nonumber \\
& = & 0.267638... \times R^2 + 0.965925... \times R.
\end{eqnarray}

Next suppose $\Gamma \in \C_2$ is WR, but not of the form $\Gamma = \Lambda(\ba)$ for some $\ba \in S(\Gamma)$, then $\Gamma = \alpha \begin{pmatrix} 1 & 1 \\ 1 & -1 \end{pmatrix} \zed^2$ for some $\alpha \in \zed_{>0}$, by Lemma~\ref{cycl-2}. Now, $\lambda_1(\Gamma) \leq R$ if and only if
$$0 < \alpha \leq \frac{R}{\sqrt{2}},$$
and so $\alpha$ can be equal to $1,2,\dots,[R/\sqrt{2}]$. Since $\alpha$ identifies $\Gamma$ uniquely,~\eqref{g-2} follows. This completes the proof.
\endproof

We are now ready to prove Theorem~\ref{cyclic_low}.

\proof[Proof of Theorem~\ref{cyclic_low}]
Notice that
$$\# \left\{ \Gamma \in \C_2: \lambda_2(\Gamma) \leq R,\ \Gamma \text{ is WR} \right\} = f_2(R)+g_2(R).$$
The result now follows directly from Lemma~\ref{cycl-2-est}.
\endproof
\bigskip

\section{Permutation invariance}
\label{perm}

Let $S_N$ be the group of permutations on $N \geq 2$ elements and define the action of $S_N$ on $\real^N$ as in~\eqref{action}. In fact, for each $\tau \in S_N$ define  $E_{\tau}$ to be the $N \times N$ matrix obtained from the $N \times N$ identity matrix $I_N$ by permuting its rows with~$\tau$; in other words, $E_{\tau} = (e_{ij})_{1 \leq i,j \leq N}$ where $e_{ij}=1$ whenever $j=\tau(i)$ and $e_{ij}=0$ otherwise. These are the well-known permutation matrices. Then for every $\bx \in \real^N$,
$$\tau \bx = E_{\tau} \bx.$$
It is easy to check that the map $\psi : S_N \to \GL_N(\zed)$ given by $\tau \mapsto E_{\tau}$ is a faithful representation of $S_N$ in $\GL_N(\real)$, and we write $\psi(S_N)$ for its image. Notice that the rotational shift operator is given precisely by the $N$-cycle $(1\ 2\dots N) \in S_N$:
\begin{equation}
\label{rot_perm}
\rot(\bx) = E_{(1\ 2 \dots N)} \bx = \begin{pmatrix} 0 & \dots & 0 & 1\\ 1 & \dots & 0 & 0\\ \vdots & \dots & \vdots & \vdots \\ 0 & \dots & 1 & 0 \end{pmatrix} \bx.
\end{equation}
Observe also that each matrix $E_{\tau}$ is orthogonal, and hence lattices $\Lambda$ and $\tau \Lambda := E_{\tau} \Lambda$ are isometric. This in particular means that $\Lambda$ is WR if and only if $\tau \Lambda$ is invariant for every~$\tau \in S_N$.

As in Section~\ref{intro}, we say that a lattice $\Lambda \subset \real^N$ is {\it $\tau$-invariant} (or invariant under~$\tau$) for a fixed $\tau \in S_N$ if $E_{\tau} \Lambda = \Lambda$. It is clear that $\Lambda$ is $\tau$-invariant if and only if it is $\sigma$-invariant for every permutation $\sigma$ in $\left< \tau \right>$, the cyclic group generated by $\tau$. This observation together with~\eqref{rot_perm} readily implies that cyclic lattices are precisely the sublattices of $\zed^N$ which are invariant under the cyclic permutation group~$\left< (1\ 2 \dots N) \right>$. Further notice that if $\Lambda$ is $\tau$-invariant and $\sigma$-invariant for some two elements $\sigma,\tau \in S_N$, then it is $(\sigma \tau)$-invariant. Recall that the transposition $(1\ 2)$ and $N$-cycle $(1\ 2\dots N)$ together generate $S_N$, and hence any cyclic lattice that is also $(1\ 2)$-invariant is invariant under the entire group~$S_N$. We can now extend our results on cyclic lattices to $\tau$-invariant full-rank sublattices of $\zed^N$ for any $N$-cycle~$\tau$.

\proof[Proof of Corollary~\ref{N-cycle}]
Let $\tau \in S_N$ be an $N$-cycle, and let us write $\sigma$ for the $N$-cycle $(1\ 2\ \dots\ N)$. Since all $N$-cycles are in the same conjugacy class, there exists $g \in S_N$ such that $\tau = g \sigma g^{-1}$. Then a lattice $\Gamma$ is $\tau$-invariant if and only if the lattice $g^{-1} \Gamma$ is $\sigma$-invariant, i.e., cyclic. Since lattices $\Gamma$ and $g^{-1} \Gamma$ are isometric, it follows that the sets
$$\left\{ \Gamma \in \C_N : \lambda_N(\Gamma) \leq R \right\},\ \left\{ \Gamma \in \C_N(\tau) : \lambda_N(\Gamma) \leq R \right\}$$
are in bijective correspondence, as are the corresponding subsets of $\WR$ and $\WR_1$ lattices, for each $R \in \real_{>0}$. The statement of the corollary now follows from Theorem~\ref{q_cyclic}.
\endproof

Since permutation invariant sublattices of $\zed^N$ are a natural generalization of cyclic lattices, we conclude with two questions about them.

\begin{quest} \label{perm_q1} Do permutation invariant full-rank sublattices of $\zed^N$ have some underlying algebraic structure? More specifically, which of them, if any, can be obtained from ideals in some polynomial rings, analogously to the construction of cyclic lattices from ideals in~$\zed[x]/(x^N -1)$?
\end{quest}

\begin{quest} \label{perm_q2} How many WR lattices are there among all $\tau$-invariant sublattices of $\zed^N$ for an arbitrary permutation $\tau \in S_N$?
\end{quest}

\noindent
A certain approach to Question~\ref{perm_q2} by means of extending the current method and studying automorphism groups of lattices is the subject of~\cite{stephan_xun}.

Both of the above questions can also be extended to {\it signed} permutation invariant lattices. Let $\J_N \cong (\zed/2\zed)^N$ be the finite abelian subgroup of $\GL_N(\zed)$ consisting of diagonal matrices with all diagonal entries being~$\pm 1$. For a fixed $g \in \J_N$ and $\tau \in S_N$, we will say that a lattice $\Lambda \subset \real^N$ is {\it $g$-signed $\tau$-invariant} if $g E_{\tau} \Lambda = \Lambda$. Now we can ask Questions~\ref{perm_q1} and~\ref{perm_q2} for signed permutation invariant lattices. As an example, let
$$g = \begin{pmatrix} -1 & 0 & \dots & 0 \\ 0 & 1 & \dots & 0 \\ \vdots & \vdots & \dots & \vdots \\ 0 & 0 & \dots & 1 \end{pmatrix} \in \J_N,\ \tau = (1\ 2\dots N) \in S_N,$$
then $g$-signed $\tau$-invariant sublattices of $\zed^N$ are images of ideals in the quotient polynomial ring $\zed[x]/(x^N +1)$ under the same map $\rho$ as for cyclic lattices in Section~\ref{intro}; we will call these the signed cyclic lattices. For instance, the signed cyclic lattices in dimension~2 are of the form
$$\begin{pmatrix} a & -b \\ b & a \end{pmatrix} \zed^2,\ a,b \in \zed.$$
These are orthogonal sublattices of $\zed^2$, which come from ideals in $\zed[x]/(x^2 +1)$ (alternatively, from ideals in Gaussian integers $\zed[i]$ under the standard Minkowski embedding of $\que(i)$ into the real plane), and are always WR. This observation suggests that signed cyclic lattices in higher dimensions may also have better than average chances of being WR.
\bigskip

{\bf Acknowledgment.} We would like to thank the referees for the highly helpful suggestions, which significantly improved the quality of the paper.
\bigskip

%\nocite{*}
\bibliographystyle{plain}  % Here the bibliography 
\bibliography{cyclic}    % is inserted.
\end{document}